\documentclass[11pt]{article}
\usepackage{amsmath,amssymb,amsbsy,amsfonts,amsthm,latexsym,amsopn,amstext,amsxtra,euscript,amscd}
\usepackage{graphicx}
\usepackage{latexsym}
\usepackage{epsfig}

\newtheorem{theorem}{Theorem}
\newtheorem{proposition}{Proposition}

\newtheorem{lemma}{Lemma}

\def\ds{\displaystyle}
\def\n{\noindent}

\smallskip
\def\bsq{\blacksquare}
\def\eproof{$\hfill\bsq$\par}

\begin{document}
\date{August 12th, 2008}
\title{\bf On independent sets in purely atomic probability spaces with geometric distribution }
\author{\sc
Eugen J. Ionascu and  Alin A. Stancu\\
\\
\small  Department of Mathematics,
\small Columbus State University\\
\small Columbus, GA 31907, USA \\ \small{\em
ionascu\_eugen@colstate.edu}, \ \ \
 \small{\em stancu\_alin1@colstate.edu}}
\maketitle
\renewcommand{\thefootnote}{}
\footnote{ \noindent\textbf{} \vskip0pt \textbf{Mathematics Subject
Classification:} Primary 60A10\vskip0pt\textbf{Key Words:}
Independence, purely atomic probability, geometric distribution}

\begin{abstract}
We are interested in constructing concrete independent events in
purely atomic probability spaces with geometric distribution. Among
other facts we prove that there are uncountable many sequences of
independent events.
\end{abstract}

\section{Introduction}

Let us assume  a fixed ratio $r$ is given, $r\in (0,1)$. In what
follows we will work with the discrete probability space $\Bbb
N_0=\{0, 1,2,3,...\}$ and the usual geometric  probability on $\cal
A$ (all subsets of $\Bbb N_0$) defined by
$$P_r(E):=\frac{1-r}{r} \sum_{k\in E\setminus \{0\} }r^{k}\ for\ every\ set\ E\in \cal A.$$
We are interested to study the class of independent sets in this
probability space. We are going to follow \cite{b} and define:
\par\vspace{0.1in}

{\it $A,B\in \Omega$ are called independent if  $P(A\cap
B)=P(A)P(B)$.}
\par\vspace{0.1in}
\n  With this definition for every $E\subset \Omega$, $\Omega$ and
$E$ are independent and $\emptyset$ and $E$ are also independent.
These are clearly trivial examples.  Three or more subsets of
$\Omega$, $A_1$, ..., $A_n$ are called {\it mutually independent} or
simply {\it independent} if for every choice of $k$ ($n\ge k\ge 2$)
such sets, say $A_{i_1}$,...,$A_{i_k}$, we have

\begin{equation}\label{nindsets}
P(\bigcap_{j=1}^k A_{i_j})=\prod_{j=1}^k P(A_{i_j}).
\end{equation}

So, for $n$ ($n\ge 2$) independent sets one needs to have $2^n-n-1$
relations as in (\ref{nindsets}) to be satisfied. An infinite family
of subsets is called independent if each finite collection of these
subsets is independent. Events are called trivial if their
probability is 0 or 1.

If $n\in \Bbb N$ then $\Omega(n)$ usually denotes the number of
primes dividing $n$ counting their multiplicities (see \cite{nzm}).
In \cite{be} and \cite{eg}, independent families of events have been
studied for finite probability spaces with uniform distribution.
Eisenberg and Ghosh \cite{eg} show that the number of nontrivial
independent events in such spaces cannot be more than $\Omega(m)$
where $m$ is the cardinality of the space. This result should be
seen in view of the known fact (see Problem 50, Section 4.1 in
\cite{gs}) that if $A_1$, $A_2$, ...., $A_n$ are independent
non-trivial events of a sample space $X$ then $|X|\ge 2^n$. One can
observe that in general $\Omega(m)$ is considerably smaller than
$\log_2 m$. It is worth mentioning that according to \cite{ess} the
first paper to deal to this problem in uniform finite probability
spaces is \cite{ss}. In their paper, Shiflett and Shultz \cite{ss}
raise the question of the existence of spaces with no non-trivial
independent pairs, called {\it dependent probability spaces}. A
space containing non-trivial independent events is called {\it
independent}. For uniform distributed probability spaces $X$, as a
result of the work in \cite{eg} and \cite{be}, $X$ is dependent if
$|X|$ a prime number and independent if $|X|$ is composite. For
denumerable sets $X$ one can see the construction given in
\cite{ess} or look at the Example 1.1 in \cite{sm}. For our spaces,
the Example 1.1 does not apply and in fact, we will construct
explicitly lots of independent sets.

For every $n\in \Bbb N$ one can consider the following space of
geometric probability distribution, denoted here by ${\cal
G}_n:=([n],{\cal P}([n]),P)$ where $P(k)=q^k$ with $k\in
[n]:=\{1,2,3,...,n\}$ and of course $q$ is the positive solution of
the equation
$$\sum_{k=1}^nq^k=1.$$ This space is independent for every $n\ge 4$
with $n$ composite. Indeed, if $n=st$ with $s,t\in \Bbb N$ $s,t\ge
2$ on can check that the sets $A:=\{1,2,3,...,s\}$,
$B:=\{1,s+1,2s+1...,(t-1)s+1\}$ represent non-trivial independent
events. To match the uniform distribution situation, it would be
interesting if ${\cal G}_n$ was a dependent space for every $n$
prime.

 The class of independent sets is important in probability
theory for various reasons. Philosophically speaking, the concept of
independence is at the heart of the axiomatic system of modern
probability theory introduced by A. N. Kolmogorov in 1933. More
recently, it was shown in \cite{crv} that two probability measures
on the same space which have the same independent (pairs of) events
must be equal if at least one of them is atomless.  This was in fact
a result of A. P. Yurachkivsky from 1989 as the same authors of
\cite{crv} point out in the addendum to their paper that appeared in
\cite{crvbis}.

On the other hand, Szekely and Mori \cite{sm} show that if the
probability space is atomic then there may be no independent sets or
one may have a sequence of such sets. The following result that
appeared in \cite{sm} is a sufficient condition for the existence of
a sequence of independent events in the probability space.

\begin{theorem}\label{primateorema}   If the range of a purely atomic probability measure
contains and interval of the form $[0,\epsilon)$ for some
$\epsilon>0$ then  there are infinitely many independent sets in the
underlying probability space.
\end{theorem}

Let us observe that, if $r=1/2$ the probability space $(\Bbb
N_0,{\cal A}, P_{1/2})$ does satisfy the hypothesis of the above
theorem with $\epsilon=1$ because every number in $[0,1]$ has a
representation in base $2$. On the other hand if, let us say
$r=1/3$, then the range of $P_{1/3}$ is the usual Cantor set which
has Lebesgue measure zero, so Theorem~\ref{primateorema} does not
apply to $(\Bbb N_0,{\cal A},P_{1/3})$. However, we will show that
there are uncountably many pairs of sets that are independent in
$(\Bbb N_0,{\cal A},P_r)$ for every $0<r<1$ (these sets do not
depend of $r$).

\section{Independent pairs of events for denumerable spaces}

The first result we would like to include is in fact a
characterization, under some restrictions of $r$, of all pairs of
independent events $(A,B)$, in which one of them, say $B$, is fixed
and of a certain form. This will show in particular that there are
uncountably many such pairs. In order to state this theorem we need
to start with a preliminary ingredient.

\begin{lemma}\label{prellemma}
For $m\ge 1$, consider the function given by
$$f(x)=(2x-1)(1+x^m)-x^m\ for\ all\ x\in [0,1].$$ The function $f$ is
strictly increasing and it has unique zero in $[0,1]$ denoted by
$t_m$. Moreover, for all $m$ we have $t_m>1/2$, the sequence
$\{t_m\}$ is decreasing and $$\ds \lim_{m\to
\infty}t_m=\frac{1}{2}.$$
\end{lemma}

Having $t_m$ defined as above we can state our first theorem.

\begin{theorem}\label{theorempairs}  For every
natural number $n\ge 2$, we define the events $E:=\{0,n-1\}$ and

\begin{equation}\label{defB}
\begin{array}{l}
B:=\{\underset{n-1}{\underbrace{1,2,...,n-1}},\underset{n-1}{\underbrace{2n-1,2n,...,3n-3}},\\
\\
\underset{n-1}{\underbrace{4n-3,4n-2,...,5n-5},...}\}.\end{array}
\end{equation}

\n Also, for $T\subset B$  an arbitrary nonempty subset we set
$A:=E+T$ with  the usual definition of addition of two sets in a
semigroup. Then $A$ and $B$ are independent events in $(\Bbb
N_0,{\cal A},P_r)$.

Conversely, if $r< t_m$ (where $m=n-1$ and $t_m$ as in
Lemma~\ref{prellemma}), $B$ is given as in (\ref{defB}) and $A$
forms an independent pair with $B$, then $A$ must be of the above
form, i.e. $A=E+T$ for some $T\subset B$.
\end{theorem}

\n {\bf Proof of Lemma~\ref{prellemma}.} The function $f$ has
derivative $f'(x)=2(1+x^m)-2m(1-x)x^{m-1}$, $x\in (0,1]$. For $m\ge
2$, using the Geometric-Arithmetic Mean inequality we have

$$(m-1)(1-x)x^{m-1}\le \left[\frac{(m-1)(1-x)+\underset{m-1}{\underbrace{x+x+...+x}}}{m}\right]^m=
\left(\frac{m-1}{m}\right)^m$$ \n and so $m(1-x)x^{m-1}\le
(\frac{m-1}{m})^{m-1}\le 1$ which implies $m(1-x)x^{m-1}\le 1$. This
last inequality is true for $m=1$ too. This implies that
$$f'(x)=2(1+x^m)-2m(1-x)x^{m-1}\ge 2x^m>0$$ for all $x\in (0,1]$.
Therefore the function $f$ is strictly increasing and because
$f(1/2)=-\frac{1}{2^m}<0$ and $f(1)=1>0$, by the Intermediate Values
Theorem there must be an unique solution $x=t_m$, of the equation
$f(x)=0$ in the interval $(1/2,1)$. Because
$f(t_{m-1})=\left(\frac{1-t_{m-1}}{1+t_{m-1}^{m-1}}\right)
t_{m-1}^{m-1}>0$ we see that $t_{m}<t_{m-1}$ for all $m\ge 2$. Since
$(2t_m-1)(1+t_m^m)=t_m^m$ we can let $m$ go to infinity in this
equality and obtain $t_m\to 1/2$.\eproof

Using Maple, we got some numerical values for the sequence $t_m$:
$t_1=\frac{1}{\sqrt{2}}\approx 0.707$, $t_2\approx 0.648$,
$t_3\approx 0.583$, $t_4\approx 0.539$ and for instance
$t_{10}\approx 0.5005$.
\vspace{0.1in}

{\bf Proof of Theorem~\ref{theorempairs}.} First let us check  that
$E_1=E+1=\{1,n\}$ and $B$ are independent. Since $E_1\cap B=\{1\}$,
$P(\{1\})=\frac{1-r}{r}r=1-r$ and
$P_r(E_1)=\frac{1-r}{r}(r+r^n)=(1-r)(1+r^{n-1})$, we have to show
that $P_r(B)=\frac{1}{1+r^{n-1}}$. We have
$$P_r(B)=\frac{1-r}{r}\left(\sum_{j=1}^m r^j\right)\left(\sum_{i=0}^{\infty}r^{2mi}\right)=
\frac{r-r^{m+1}}{r}\frac{1}{1-r^{2m}}=\frac{1}{1+r^{m}}$$ \n which
is what we needed. Now, suppose $b\in B$ and consider
$E_b=E+b=\{b,b+n-1\}$. We notice that by the definition of $B$, the
intersection $B\cap E_b$ is $\{b\}$. Hence, $P_r(B\cap
E_b)=\frac{1-r}{r}r^b=(1-r)r^c$ (with $c=b-1$) and
$$P_r(B)P_r(E_b)=\frac{1}{1+r^m}\frac{1-r}{r}\left(r^b+r^{b+m}\right)=(1-r)r^c.$$

\n Hence, $B$ and $E_b$ are independent for every $b\in B$.

\begin{quote}
{\it Next we would like to observe that if $(F_1,B)$ and $(F_2,B)$
are independent pairs of events and $F_1\cap F_2=\emptyset$, then
$F_1\cup F_2$ and  $B$ are independent events as well.}
\end{quote}
 Indeed, by the  given assumption we can write
$$
\begin{array}{l}
P_r(B\cap (F_1\cup F_2))=P_r((B\cap F_1)\cup (B\cap F_2))=P_r(B\cap
F_1)+P_r(B\cap F_2)=\\ \\
P_r(B)P_r(F_1)+P_r(B)P_r(F_2)=P_r(B)(P_r(F_1)+P_r(F_2))=P_r(B)P_r(F_1\cup
F_2).
\end{array}
$$

In fact, the above statement can be generalized to a sequence of
sets $F_k$ which are pairwise disjoint, due to the fact that $P_r$
is a genuine finite measure and so it is continuous (from below and
above). Then if $T\subset B$ is nonempty, $A=E+T=\bigcup_{b\in
B}E_b$ is countable union and since $E_b\cap E_{b'}=\emptyset$ for
all $b,b'\in B$ ($b\not=b'$) the above observation can be applied to
$\{E_b\}_{b\in T}$. So, we get that $B$ and $A$ are independent.

For the converse, we need the following lemma.

 \begin{lemma}\label{adoualema} If $L \subset \Bbb N_0 \setminus B$ and the smallest element of
 $L$ is $s=(2i-1)m+j$, where $i,j \in \Bbb N$, $j\le m$,  then $$P_r(L)\le
r^{s-1}-\frac{r^{2im}}{1+r^m}.$$
 \end{lemma}

\n {\bf Proof of Lemma~\ref{adoualema}} Indeed, we have

$$
\begin{array}{c}
P_r(L)\le \frac{1-r}{r}[(r^s+r^{s+1}+...+r^{2im})+(r^{(2i+1)m+1}+....)]=\\
\\
r^{s-1}-r^{2im}+r^{2im}P_r(\Omega\setminus
B)=r^{s-1}-r^{2im}+r^{2im}(1-\frac{1}{1+r^m})=r^{s-1}-\frac{r^{2im}}{1+r^m}.
\end{array}
$$
\eproof

So, let us assume that $r<t_m$, $B$ is as in (\ref{defB}) and $A$ is
independent of $B$. We let $T$ be the intersection of $A$ and $B$
and we put $\alpha:=P_r(T)/P_r(B)$. Also, define $A':=T+\{0,n-1\}$,
$L=A\setminus A'$ and $L'=A'\setminus A$. We have clearly
$L,L'\subset \Omega\setminus B$. By the first part of our theorem
$P_r(A')=\alpha$. Because $A$ and $B$ are independent $P_r(A)$ must
be equal to $\alpha$ as well. Hence $P_r(A)=P_r(A')$ which attracts

\begin{equation}\label{important}\sum_{k\in L'}r^k=\sum_{k\in L}r^k \Leftrightarrow \sum_{k\in
L\cup L'}r^k=2\sum_{k\in L' }r^k.
\end{equation}

From (\ref{important}), it is clear that $L'=\emptyset$ if an only
if $L=\emptyset$ and so if $L'$ is empty then $A=A'$, which is what
we need in order to conclude our proof. By way of contradiction,
suppose $L'\not=\emptyset$ (or equivalently $L\not=\emptyset$) we
can assume without loss of generality that $L'$ contains the
smallest number of $L'\cup L$, say $s$ which is written as in
Lemma~\ref{adoualema}. Thus from equality (\ref{important}) we have
$P_r(L\cup L')\ge 2P_r(L')$ and then by Lemma~\ref{adoualema} we get

$$r^{s-1}-\frac{r^{2im}}{1+r^m}\ge 2(1-r)r^{s-1}
\Leftrightarrow 2r\ge 1+  \frac{r^{2im+1-s}}{1+r^m} \Leftrightarrow
2r\ge 1+\frac{r^{n-j}}{1+r^m}.$$

Therefore for every $n$ and $1\le j\le m$, $$2r\ge
1+\frac{r^{n-j}}{1+r^m}\ge 1+\frac{r^{m}}{1+r^m} \Rightarrow
f(r)=(2r-1)(1+r^m)-r^m\ge 0.$$

By Lemma~\ref{prellemma} we see that $r\ge t_m$ which is a
contradiction. It remains that $L$ and $L'$ must be empty and so
$A=A'$. \eproof

In the previous theorem, since $T$ was an arbitrary subset of an
infinite set we obtain an uncountable family of pairs of independent
sets.

\vspace{0.1in}

{\bf Remark 1:} If $r=\sqrt{\frac{1}{\phi}}$ where $\phi$ stands for
the classical notation of the golden ratio (i.e.
$\phi=\frac{\sqrt{5}+1}{2}$), $n=2$, $B=\{1,3,5,7,...\}$ as in
(\ref{defB}), and $A=\{1,4,6\}$, then one can check that
$P_r(B)=\frac{1}{1+r}$, $P_r(A\cap B)=1-r$,
$P_r(A)=(1-r)(1+r^3+r^5)$. So the equality $P_r(A\cap
B)=P_r(A)P_r(B)$ is equivalent to $1+r=1+r^3+r^5$ which is the same
as $r^4+r^2-1=0$. One can easily see that this last equation is
satisfied by $r=\sqrt{\frac{1}{\phi}}$. Hence $A$ and $B$ are
independent but clearly $A$ is not a translation of $\{0,1\}$ with a
subset of $B$. Therefore the converse part in
Theorem~\ref{primateorema} cannot  be extended to numbers $r\ge t_m$
such as $r=\sqrt{\frac{1}{\phi}}$. In fact, we believe that the
constants $t_m$ are sharp, in the sense that for all $r>t_m$ the
converse part is false,  but an argument for showing this is beyond
the scope of this paper.

\vspace{0.1in}

{\bf Remark 2:} Another family of independent events which seems to
have no connection with the ones constructed so far is given by
$A=\{1,2,3,4,...,n-1,n\}$ and $B=\{n,2n,3n,...\}$, with $n\in \Bbb
N$. A natural question arises as a result of this wealth of
independent events: can one characterize all pairs $(A,B)$ which are
independent regardless the value of the parameter $r$?

\section{Three independent events}

The next theorem deals with the situation in which two sets as in
the construction of Theorem~\ref{theorempairs} form with $B$ given
by (\ref{defB}), a triple of independent sets.

Let us observe that if $A_1$, $A_2$, and $B$ are mutually
independent then by Theorem~\ref{theorempairs} (at least if $r\in
(0,t_m)$), $A_1$ and $A_2$ must be given by $A_i=T_i+E$ with
$T_i\subset B$, $i=1,2$. Therefore $A_1\cap A_2=(T_1\cap T_2)+E$.

Also, we note that $P_r(A_i)=P_r(T_i)(1+r^{n-1})$, $i=1,2$, and
$P_r(A_1\cap A_2)=P_r(T_1\cap T_2)(1+r^{n-1})$. This means that the
equality $P_r(A_1\cap A_2)=P_r(A_1)P(A_2)$ is equivalent to

\begin{equation}\label{newindepdef}
P_r(T_1\cap T_2)=P_r(T_1)P_r(T_2)(1+r^{n-1}).
\end{equation}

On the other hand the condition $P_r(A_1\cap A_2\cap
B)=P_r(A_1)P_r(A_2)P(B)$ reduces to $$P_r(T_1\cap
T_2)=P_r(T_1)P_r(T_2)(1+r^{n-1})^2P_r(B),$$ which is the same as
(\ref{newindepdef}). So, three sets $A_1$, $A_2$ and $B$ are
independent if and only if (\ref{newindepdef}) is satisfied.  Let us
notice that the condition (\ref{newindepdef}) may be interpreted as
a conditional probability independence relation:

\begin{equation}\label{newindepascondpr}
P_r(T_1\cap T_2|B)=P_r(T_1|B)P_r(T_2|B).
\end{equation}

At this point the construction we have in Theorem~\ref{theorempairs}
can be repeated. As a result, regardless of what $r$ is, we obtain
an uncountable family of there events which are mutually independent
in $(\Bbb N_0, {\cal A}, P_r)$.

\begin{theorem} For a fixed $n\ge 3$, we consider $B$ as in (\ref{defB}), and pick $b\in \{2,...,n-1\}$ such that
$2(b-1)$ divides $m=n-1$ ($m=2(b-1)k$). For  $F:=\{0,b-1\}$, we let

\begin{equation}\label{defbunu}
\begin{array}{c}
B_1':=\{\underset{b-1}{\underbrace{1,2,...,b-1}},\underset{b-1}{\underbrace{2b-1,2b,...,3b-3}},
\underset{b-1}{\underbrace{4b-3,4b-2,...,5b-5}},\\
\\...,\underset{b-1}{\underbrace{(2k-2)(b-1)+1,...,(2k-1)(b-1)}}\},\end{array}
\end{equation}

$$B_1:=B_1'\cup (B_1'+2m)\cup (B_1'+4m)\cup (B_1'+6m)\cup ...$$

\n and $T$ a subset of $B_1$. Then $T_1:=F+T$ and $B_1$ are
independent sets relative to the induced probability measure on $B$.
Moreover, $A_1:=T_1+\{0,n-1\}$, $A_2:=B_1+\{0,n-1\}$ and $B$ form a
triple of mutually independent sets in $(\Bbb N_0,{\cal A},P_r)$ for
all $r$.
\end{theorem}

\proof The second part of the theorem follows from the
considerations we made before the theorem and from the first part.
To show the first part we need to check (\ref{newindepdef}) for
$T_1$ and $T_2=B_1$. Let us remember that
$$
\begin{array}{l}
B=\{\underset{n-1}{\underbrace{1,2,...,n-1}},\underset{n-1}{\underbrace{2n-1,2n,...,3n-3}},\\
\\
\underset{n-1}{\underbrace{4n-3,4n-2,...,5n-5},...}\},\ and \ \
P_r(B)=\frac{1}{1+r^m}.\end{array}
$$

We observe that $B_1'\subset \{1,2,...,n-1\}$ and so $B_1\subset B$.
Let us first take into consideration the case $T=\{1\}$. Since
$T_1=\{1,b\}$ we get $T_1\cap T_2=\{1 \}$,
$P_r(T_1)=(1-r)(1+r^{b-1})$, and
$$P_r(B_1)=P_r(B_1')(1+r^{2m}+r^{4m}+r^{6m}+...)=\frac{P_r(B_1')}{1-r^{2m}}.$$

So, it remains to calculate $P_r(B_1')$:

$$
\begin{array}{l}
P_r(B_1')=\frac{1-r}{r}(r+r^2+...r^{b-1})(1+r^{2(b-1)}+r^{4(b-1)}+...+r^{2(k-1)(b-1)})\\
\\

\ds
=(1-r^{b-1})\frac{1-r^{2k(b-1)}}{1-r^{2(b-1)}}=\frac{1-r^{m}}{1+r^{b-1}}\Rightarrow
P_r(B_1)=\frac{1}{(1+r^{b-1})(1+r^{m}).}
\end{array}
$$

This shows that (\ref{newindepdef}) is satisfied. In the general
case, i.e. $T$ an arbitrary subset of $B_1$, we proceed as in the
proof of Theorem~\ref{primateorema}. \eproof

\section{Uncountable sequences of independent events}
In \cite{sm}, Szekely and Mori give an example of an infinite
sequence of independent sets in $(\Bbb N_0,{\cal A},P_{1/2})$. Given
an infinite sequence of independent sets $\{A_n\}_n$ we may assume
that $P_r(A_k)\le \frac{1}{2}$ and so by Proposition 1.1 in
\cite{sm} we must have

$$\sum_{k=1}^{\infty}P_r(A_k)<\infty.$$

Let us observe that Theorem~\ref{theorempairs} can be applied to a
different space now that can be constructed  within $B$ given by
(\ref{defB})in terms of classes: $\widehat{\Bbb N}_0=\{
\hat{0},\hat{1},\hat{2},...\}$ where $\hat{0}=\emptyset$,
$\hat{1}:=\{1,2,...,n-1\}$, $\hat{2}:=\{2n-1,2n,...,3n-3\}$,
$\hat{3}:=\{4n-3,4n-2,...,5n-5\}$, ..., and the probability on this
space is the conditional probability as subsets of $B$.

Hence for $k\in \Bbb N$, one can check that
$$P(\hat{k})=\frac{1-r^{2m}}{r^{2m}}r^{2km},\ with \ m=n-1.$$ This
shows that this space is isomorphic to $(\Bbb N_0,{\cal A},P_s)$
with $s=r^{2m}$.

One can check by induction the following proposition.

\begin{proposition}\label{induction} Let $n\in \Bbb N$, $n\ge 2$. If $A_1$,...,$A_n$ are
independent in $\widehat{\Bbb N}_0$ then $A_1+T$, $A_2+T$,...,
$A_n+T$ and $B$ are indepenedent in $(\Bbb N_0,{\cal A},P_r)$.
\end{proposition}

This construction can be then iterated indefinitely giving rise of a
sequence $B$, $B_1$, $B_2$,...,  which is going to be independent
and its construction is in terms of a sequence $(n,n_1,n_2,...)$
with $n_k\ge 2$. As a result, we have a countable way of
constructing sequences of independent sets. This construction
coincides with the one in \cite{sm} if $n_k=2$ for all $k\in \Bbb
N$.

\end{document}